\documentclass[twoside,12pt]{article}
\usepackage{amsmath}
\usepackage{amssymb}

\def\q{\quad}
\def\qq{\qquad}

\def\mod#1{\ (\text{\rm mod}\ #1)}
\def\t{\text}
\def\f{\frac}
\def\e{\equiv}
\def\b{\binom}

\def\sls#1#2{(\f{#1}{#2})}
 \def\ls#1#2{\big(\f{#1}{#2}\big)}
\def\Ls#1#2{\Big(\f{#1}{#2}\Big)}
\def\xp{\langle x\rangle_p}
\def\v2{\vskip0.2cm}

\allowdisplaybreaks

\begin{document}

 \centerline {\large\bf
Curious identities involving Legendre polynomials and Ap\'ery-like numbers}
\par\q\newline
\centerline{Li-Li Cui$^1$ and Zhi-Hong Sun$^2$}
\vskip0.2cm
\centerline{$^1$School of Mathematics and Statistics}
\centerline{Huaiyin Normal University}
\centerline{Huaian, Jiangsu 223300, P.R. China}
 \centerline{Email:1655250018 @stu.hytc.edu.cn}
 \vskip0.2cm
\centerline{$^2$School of Mathematics and Statistics}
\centerline{Huaiyin Normal University}
\centerline{Huaian, Jiangsu 223300, P.R. China}
 \centerline{Email: zhsun@hytc.edu.cn}
 \centerline{URL: https://maths.hytc.edu.cn/szh.htm}

\par\q
\par {\bf Abstract.} In this paper, we establish some curious identities involving Legendre polynomials and the first kind of Ap\'ery-like numbers. As applications, many new supercongruences are deduced.
\par\q

\section*{1. Introduction}
\par The famous Legendre polynomials $\{P_n(x)\}$ are given by
$$P_0(x)=1,\  P_1(x)=x\ \t{and}\  (n+1)P_{n+1}(x)=(2n+1)xP_n(x)-nP_{n-1}(x)\ (n\ge 1).\eqno{(1.1)}$$
It is well known (see for example [17]) that
$$ P_n(x)=\sum_{k=0}^n\b nk\b{n+k}k\Ls{x-1}2^k=\f 1{2^n}\sum_{k=0}^{[n/2]}
 \b nk\b{2n-2k}n(-1)^kx^{n-2k},$$
 where $[t]$ is the greatest integer not exceeding $t$. Hence,
 $P_n(-x)=(-1)^nP_n(x)$ and
 $$P_n(1)=1,\ P_n(-1)=(-1)^n,\ P_n(0)=\begin{cases} 0,&\t{if $2\nmid n$,}
 \\\f{(-1)^{n/2}}{2^n}\b n{n/2},&\t{if $2\mid n$.}
 \end{cases}\eqno{(1.2)}$$
Let $p$ be a positive integer and $x\not=y$. The famous Christofffel-Darboux formula (see [1, Remark 5.2.2] and [2, Theorem 2.3.3]) states that
$$\sum_{n=0}^{p-1}(2n+1)P_n(x)P_n(y)=\f p{x-y}
(P_p(x)P_{p-1}(y)-P_{p-1}(x)P_p(y)).\eqno{(1.3)}$$
In Section 2, we establish closed formulas for
\begin{align*}
&\sum_{n=0}^{p-1}(2n+1)^3P_n(x)P_n(y),\ \sum_{n=0}^{p-1}(2n+1)^3P_n(x)^2, \ \sum_{n=0}^{p-1}(2n+1)^5P_n(x)P_n(y),
 \\&\sum_{n=0}^{p-1}(2n+1)^5P_n(x)^2,\ \sum_{n=0}^{p-1}(2n+1)^7P_n(x),\  \sum_{n=0}^{p-1}(2n+1)^9P_n(x).\end{align*}

\par\q Let $\Bbb Z,\Bbb Z^+$ and $\Bbb R$ be the set of integers, the set of positive integers and the set of real numbers, respectively. For $a,b,c\in\Bbb R$ with $c\not=0$ let $\{u_n(a,b,c)\}$ be given by
\begin{align*}u_0=1, \ u_1=\ b,\ (n+1)^3u_{n+1}
=(2n+1)(an(n+1)+b)u_n-cn^3u_{n-1}\  (n\ge 1).\tag {1.4}\end{align*}
 If $a,b,c\in\Bbb Z$ and $u_n\in\Bbb Z$ for $n=1,2,3,\ldots$, following [16] we call $\{u_n\}$ the first kind of Ap\'ery-like numbers. Let
\begin{align*}
 &A_n= \sum_{k=0}^n\binom nk^2\binom{n+k}k^2,
\q D_n=\sum_{k=0}^n\b nk^2\b{2k}k\b{2n-2k}{n-k},\\&
b_n=\sum_{k=0}^{[n/3]}\b{2k}k\b{3k}k\b n{3k}\b{n+k}k(-3)^{n-3k}, \
T_n=\sum_{k=0}^n\b nk^2\b{2k}n^2, \
\\\qq\qq& V_n=\sum_{k=0}^n\b
nk\b{n+k}k(-1)^k\b{2k}k^216^{n-k}
 =\sum_{k=0}^n\b{2k}k^2\b{2n-2k}{n-k}^2.\end{align*}
Then $\{A_n\}$, $\{D_n\}$, $\{b_n\}$, $\{T_n\}$ and $\{V_n\}$ are
 the first kind of Ap\'ery-like numbers with
$(a,b,c)=(17,5,1),(10,4,64),(-7,-3,81),(12,4,16)$ and $(16,8,256)$,
respectively. The numbers $\{A_n\}$, $\{D_n\}$ and $\{b_n\}$ are
called Ap\'ery numbers, Domb numbers and Almkvist-Zudilin numbers,
respectively. For $\{A_n\}$, $\{D_n\}$, $\{b_n\}$, $\{T_n\}$ and
$\{V_n\}$ see A005259, A002895, A125143, A290575 and A036917 in [15].
\par Let
$$V_n(x)=\sum_{k=0}^n\b nk\b{n+k}k(-1)^k\b
xk\b{-1-x}k\ (n=0,1,2,\ldots).$$
Then $V_n=16^nV_n(-\f 12)$. From [17,(22.10) and Theorem 22.44],
$$V_n(x)=\sum_{k=0}^n\b xk^2\b{-1-x}{n-k}^2=u_n(1,2x^2+2x+1,1).$$

\par Let $p\in\Bbb Z^+$ and $a,b,d\in\Bbb R$ with $d\not=0$. In Section 3, for $r=3,5,7$ we reveal the connections between
$\sum_{n=0}^{p-1}(2n+1)^r\frac{u_n(a,b,d^2)}{d^n}$ and
$\sum_{n=0}^{p-1}(2n+1)\frac{u_n(a,b,d^2)}{d^n}$.
\par Let $p$ be a prime greater than $3$. In Section 4, using the identities in Section 3 we deduce the congruences for
\begin{align*}&\sum_{n=0}^{p-1}(2n+1)^5A_n,\ \sum_{n=0}^{p-1}(2n+1)^7A_n,\ \sum_{n=0}^{p-1}(2n+1)^5\f{T_n}{4^n},\
\\&\sum_{n=0}^{p-1}(2n+1)^7\f{T_n}{4^n},\
\sum_{n=0}^{p-1}(2n+1)^5\f{D_n}{8^n},\ \sum_{n=0}^{p-1}(2n+1)^7\f{D_n}{8^n}
\end{align*} modulo $p^5$, and the congruence for
\begin{align*}&\sum_{n=0}^{p-1}(2n+1)^5\f{b_n}{(-9)^n},\ \sum_{n=0}^{p-1}(2n+1)^7\f{b_n}{(-9)^n}\pmod{p^3},
\\&\sum_{n=0}^{p-1}(2n+1)^5\f{T_n}{(-4)^n},\ \sum_{n=0}^{p-1}(2n+1)^7\f{T_n}{(-4)^n}\pmod{p^4}.
\end{align*}
We also obtain the congruences for
$\sum_{n=0}^{p-1}(2n+1)^5(-1)^nV_n(x)$ and
$\sum_{n=0}^{p-1}(2n+1)^7(-1)^nV_n(x)$ modulo $p^4$.

\par In addition to the above notation, throughout this paper we also use the following notations. For an odd prime $p$ let $\Bbb Z_p$ be the set of those rational numbers whose denominators are not divisible by $p$.
 For $a\in\Bbb Z$ and given odd prime $p$, $\ls ap$ denotes the Legendre symbol. The Bernoulli numbers $\{B_n\}$, Euler numbers $\{E_n\}$ and the sequence
$\{U_n\}$ are defined by
\begin{align*} &B_0=1,\q\sum_{k=0}^{n-1}\b nkB_k=0\q(n\ge 2),\\& E_{2n-1}=0,\q E_0=1,\q E_{2n}=-\sum_{k=0}^{n-1}
\b {2n}{2k}E_{2k}\q(n\ge 1),
\\& U_{2n-1}=0,\q U_0=1,\q U_{2n}=-2\sum_{k=0}^{n-1}
\b {2n}{2k}U_{2k}\q(n\ge 1).
\end{align*}
The Euler
polynomials $\{E_n(x)\}$ are given by
$$E_n(x)=\f 1{2^n}\sum_{k=0}^n\b nk(2x-1)^{n-k}E_k.$$

\section*{\bf 2. Identities for Legendre polynomials}
\par\q Putting $p_n(x)=P_n(x)/\sqrt{\f 2{2n+1}}$ in [2, Corollary 2.3.4]
and then applying the known formula $(1-x^2)P'_n(x)=nP_{n-1}(x)-nxP_n(x)$ yields the formula
$$\sum_{n=0}^{p-1}(2n+1)P_n(x)^2
=\f{p^2}{1-x^2}\big(P_p(x)^2-2xP_{p-1}(x)P_p(x)+P_{p-1}(x)^2\big).
\eqno{(2.1)}$$
Putting $r=1,\ c=1$ and $b(n)=(2n+1)x$ in [16, Theorem 3.3] gives
$$\sum_{n=0}^{p-1}(2n+1)(-1)^nP_n(x)^2=(-1)^{p-1}\f pxP_p(x)P_{p-1}(x).\eqno{(2.2)}$$
Combining (2.1) and (2.2) yields
\begin{align*}&\sum_{n=0}^{[(p-1)/2]}(4n+1)P_{2n}(x)^2
\\&=\f{p^2}{2(1-x^2)}
(P_p(x)^2+P_{p-1}(x)^2)+\Big((-1)^{p-1}\f p{2x}-\f{p^2x}{1-x^2}\Big)P_p(x)P_{p-1}(x).\tag{2.3}
\end{align*}
\par Inspired by (1.3) and (2.1), with the help of Maple we find the following identities.
\v2
\par{\bf Theorem 2.1} {\sl Let $p$ be a positive integer.
\par $(\t{\rm i})$ For $x\not=\pm 1$ we have}
\begin{align*}&\sum_{n=0}^{p-1}(2n+1)^3P_n(x)^2
\\&=\f{p^2}{3(1-x^2)^2}
\Big(((4p(p-1)+3)(1-x^2)-8)P_p(x)^2
\\&\q+((4p(p+1)+3)(1-x^2)-8)P_{p-1}(x)^2+2x((1-4p^2)(1-x^2)+8)
P_p(x)P_{p-1}(x)\Big).
\end{align*}
\par $(\t{\rm ii})$ For $x\not=y$ we have
$$\sum_{n=0}^{p-1}(2n+1)^3P_n(x)P_n(y)=g_p(x,y),$$
where
\begin{align*}g_p(x,y)&=\f p{(x-y)^3}\Big(
\big((2p+1)^2(x-y)^2-8(p+1)x(x-y)-8(1-x^2)\big)
P_p(x)P_{p-1}(y)
\\&\q-\big((2p-1)^2(x-y)^2+8(p-1)x(x-y)-8(1-x^2)\big)
P_{p-1}(x)P_p(y)
\\&\q+8p(x-y)(P_p(x)P_p(y)+P_{p-1}(x)P_{p-1}(y))\Big).
\end{align*}

\par{\it Proof.} Set
\begin{align*}
h_n(x)&=\f{n^2}{3(1-x^2)^2}
\big(((4n(n-1)+3)(1-x^2)-8)P_n(x)^2
\\&\q+((4n(n+1)+3)(1-x^2)-8)P_{n-1}(x)^2
\\&\q+2x((1-4n^2)(1-x^2)+8)
P_n(x)P_{n-1}(x)\big).
\end{align*}
Since $(n+1)P_{n+1}(x)=(2n+1)xP_n(x)-nP_{n-1}(x)$, we see that
\begin{align*}
h_{n+1}(x)&=\f{(n+1)^2}{3(1-x^2)^2}
\Big(((4(n+1)n+3)(1-x^2)-8)P_{n+1}(x)^2
\\&\q+((4(n+1)(n+2)+3)(1-x^2)-8)P_n(x)^2
\\&\q+2x((1-4(n+1)^2)(1-x^2)+8)
P_{n+1}(x)P_n(x)\Big)
\\&=\f{1}{3(1-x^2)^2}
\Big(((4(n+1)n+3)(1-x^2)-8)((2n+1)xP_n(x)-nP_{n-1}(x))^2
\\&\q+(n+1)^2((4(n+1)(n+2)+3)(1-x^2)-8)P_n(x)^2
\\&\q+2x((1-4(n+1)^2)(1-x^2)+8)
(n+1)((2n+1)xP_n(x)-nP_{n-1}(x))P_n(x)\Big)
\\&=h_n(x)+(2n+1)^3P_n(x)^2.
\end{align*}
Thus,
$$\sum_{n=0}^{p-1}(2n+1)^3P_n(x)^2=\sum_{n=0}^{p-1}(h_{n+1}(x)-h_n(x))
=h_p(x)-h_0(x)=h_p(x),$$
which proves (i).
\par Now we consider (ii). Since $(n+1)P_{n+1}(x)=(2n+1)xP_n(x)-nP_{n-1}(x)$, we see that
\begin{align*}
&g_{n+1}(x,y)\\&=\f {n+1}{(x-y)^2}\Big(
\Big((2n+3)^2(x-y)-8(n+2)x-\f{8(1-x^2)}{x-y}\Big)
P_{n+1}(x)P_{n}(y)
\\&\q-\Big((2n+1)^2(x-y)+8nx-\f{8(1-x^2)}{x-y}\Big)P_n(x)P_{n+1}(y)
\\&\q+8(n+1)(P_{n+1}(x)P_{n+1}(y)+P_n(x)P_n(y))\Big)
\\&=\f 1{(x-y)^2}\Big(
\Big((2n+3)^2(x-y)-8(n+2)x-\f{8(1-x^2)}{x-y}\Big)
((2n+1)xP_n(x)-nP_{n-1}(x))P_{n}(y)
\\&\q-\Big((2n+1)^2(x-y)+8nx-\f{8(1-x^2)}{x-y}\Big)P_n(x)
((2n+1)yP_n(y)-nP_{n-1}(y))
\\&\q+8(((2n+1)xP_n(x)-nP_{n-1}(x))((2n+1)yP_n(y)-nP_{n-1}(y))
+(n+1)^2P_n(x)P_n(y))\Big)
\\&=g_n(x,y)+(2n+1)^3P_n(x)P_n(y).
\end{align*}
Thus,
$$\sum_{n=0}^{p-1}(2n+1)^3P_n(x)P_n(y)=\sum_{n=0}^{p-1}(g_{n+1}(x,y)
-g_n(x,y))=g_p(x,y)-g_0(x,y)=g_p(x,y).$$
This proves (ii). The proof is now complete.
 \vskip0.2cm
\par{\bf Corollary 2.1} {\sl Let $p\in\Bbb Z^+$ and $x\not=1$. Then}
\begin{align*}&\sum_{n=0}^{p-1}(2n+1)^3P_n(x)
\\&=\f p{(1-x)^2}\big(((2p+1)^2(1-x)-8)P_{p-1}(x)
-((2p-1)^2(1-x)-8)P_p(x)\big)
\tag{2.4}\end{align*}
and for odd $p$ and $x\not=\pm 1$,
\begin{align*}
&\sum_{n=0}^{(p-1)/2}(4n+1)^3P_{2n}(x)
\\&=\Big(\f{p(2p+1)^2}{1-x^2}-\f{8p(1+x^2)}{(1-x^2)^2}\Big)
P_{p-1}(x)
+\Big(-\f{p(2p-1)^2x}{1-x^2}+\f{16px}{(1-x^2)^2}\Big)P_p(x).
\tag {2.5}\end{align*}
 \vskip0.2cm
 \par{\it Proof.} Putting $y=1$ in Theorem 2.1(ii) yields
 (2.4). Since $P_n(-x)=(-1)^nP_n(x)$, from (2.4) we see that
\begin{align*}&(1+x)^2\sum_{n=0}^{p-1}(2n+1)^3(-1)^nP_n(x)
\\&=(-1)^{p-1}p
\big((2p+1)^2(1+x)-8)P_{p-1}(x)+((2p-1)^2(1+x)-8)P_p(x)\big).
\end{align*}
Thus, for odd $p$ and $x\not=\pm 1$,
\begin{align*}
&2\sum_{n=0}^{(p-1)/2}(4n+1)^3P_{2n}(x)
\\&=\sum_{n=0}^{p-1}(2n+1)^3(1+(-1)^n)P_n(x)
\\&=\f{p\big((2p+1)^2(1-x)-8)P_{p-1}(x)-((2p-1)^2(1-x)-8)P_p(x)}
{(1-x)^2}
\\&\q+\f{p
\big((2p+1)^2(1+x)-8)P_{p-1}(x)+((2p-1)^2(1+x)-8)P_p(x)\big)}
{(1+x)^2}
\\&=2\Big(\f{p(2p+1)^2}{1-x^2}-\f{8p(1+x^2)}{(1-x^2)^2}\Big)
P_{p-1}(x)
+2\Big(-\f{p(2p-1)^2x}{1-x^2}+\f{16px}{(1-x^2)^2}\Big)P_p(x),
\end{align*}
which proves (2.5). The proof is now complete.
\vskip0.2cm
\par{\bf Corollary 2.2} {\sl For $p=1,3,5,\ldots$ we have}
$$\sum_{n=0}^{(p-1)/2}(4n+1)^3\f{\b{2n}n}{(-4)^n}=p(4p^2+4p-7)
\f{\b{p-1}{(p-1)/2}}{(-4)^{(p-1)/2}}.$$
Hence, for any prime $p>3$ we have
$$\sum_{n=0}^{(p-1)/2}(4n+1)^3\f{\b{2n}n}{(-4)^n}\e
(4p^3+4p^2-7p)2^{p-1}-\f 7{12}p^4B_{p-3}\mod {p^5}.$$
\par{\it Proof.} Putting $t=0$ in Corollary 2.1 and noting that
$P_{2n}(0)=\b{2n}n/(-4)^n$ and
$$(-1)^{\f{p-1}2}\b{p-1}{(p-1)/2}\e 4^{p-1}+\f{p^3}{12}B_{p-3}\mod{p^4}\eqno{(2.6)}$$
(see [17,pp.330-331])
gives the result.
\vskip0.2cm
\par{\bf Corollary 2.3} {\sl Let $p\in\Bbb Z^+$ and $x\not=0$. Then}
\begin{align*}&\sum_{n=0}^{p-1}(2n+1)^3(-1)^nP_n(x)^2
\\&=(-1)^{p-1}\f p{x^3}\big(2pxP_{p-1}(x)^2+((4p^2-1)x^2-2)P_p(x)P_{p-1}(x)-2pxP_p(x)^2
\big)\tag{2.7}\end{align*}
and for $x\not=\pm 1$,
\begin{align*}&\sum_{n=0}^{[(p-1)/2]}(4n+1)^3P_{2n}(x)^2
\\&=\Big(\f{p^2}{6(1-x^2)^2}((4p(p-1)+3)(1-x^2)-8)-(-1)^{p-1}\f{p^2}{x^2}\Big)P_p(x)^2
\\&\q+\Big(\f{p^2x}{3(1-x^2)^2}((1-4p^2)(1-x^2)+8)+(-1)^{p-1}
\f p{2x^3}((4p^2-1)x^2-2)\Big)P_p(x)P_{p-1}(x)
\\&\q+\Big(\f{p^2}{6(1-x^2)^2}((4p(p+1)+3)(1-x^2)-8)+(-1)^{p-1}\f{p^2}{x^2}
\Big)P_{p-1}(x)^2.
\end{align*}
\par{\it Proof.} Taking $y=-x$ in Theorem 2.1(ii) and noting that $P_n(-x)=(-1)^nP_n(x)$ yields (2.7). Since
 $$2\sum_{n=0}^{[(p-1)/2]}(4n+1)^3P_{2n}(x)^2=\sum_{n=0}^{p-1}
 (2n+1)^3(1+(-1)^n)P_n(x)^2,$$
 combining (2.7) with Theorem 2.1(i) yields the remaining result.
\v2
\par{\bf Theorem 2.2} {\sl Let $p\in\Bbb Z^+$ and $x\not=y$. Then}
\begin{align*}
&(x-y)^5\sum_{n=0}^{p-1}(2n+1)^5P_n(x)P_n(y)
\\&=16p^2\big((4p(p-1)+5)(x-y)^3-24x(x-y)^2+24(x^2-1)(x-y)\big)P_p(x)P_p(y)
\\&\q+16p^2\big((4p(p+1)+5)(x-y)^3-24x(x-y)^2+24(x^2-1)(x-y)\big)
P_{p-1}(x)P_{p-1}(y)
\\&\q+p\big((2p+1)^4(x-y)^4-16(p+1)(4p^2+8p+5)x(x-y)^3
\\&\q+16((12p^2+36p+29)x^2-(12p^2+12p+13))(x-y)^2
\\&\q-384(p+2)(x^3-x)(x-y)+384(x^2-1)^2\big)P_p(x)P_{p-1}(y)
\\&\q-p\big((2p-1)^4(x-y)^4+16(p-1)(4p^2-8p+5)x(x-y)^3
\\&\q+16((12p^2-36p+29)x^2-(12p^2-12p+13))(x-y)^2
\\&\q+384(p-2)(x^3-x)(x-y)+384(x^2-1)^2\big)P_{p-1}(x)P_{p}(y).
\end{align*}
\par{\it Proof.} The proof is similar to the proof of Theorem 2.1(ii). Let $f_p(x,y)$ be the right-hand side of the identity in Theorem 2.2. Using the relation $(n+1)P_{n+1}(x)=(2n+1)xP_n(x)-nP_{n-1}(x)$ one may deduce that
$$f_{n+1}(x,y)-f_n(x,y)=(x-y)^5(2n+1)^5P_n(x)P_n(y).$$
Thus,
\begin{align*}&(x-y)^5\sum_{n=0}^{p-1}(2n+1)^5P_n(x)P_n(y)
\\&=\sum_{n=0}^{p-1}(f_{n+1}(x,y)
-f_n(x,y))=f_p(x,y)-f_0(x,y)=f_P(x,y).\end{align*}
This proves the theorem.

\v2
\par{\bf Corollary 2.4} {\sl Let $p\in\Bbb Z^+$ and $x\not=0$. Then}
\begin{align*}\sum_{n=0}^{p-1}(2n+1)^5(-1)^nP_n(x)^2
&=(-1)^{p-1}\f p{x^5}\Big(((-16p^3+16p^2+4p)x^3+24px)P_p(x)^2
\\&\q+((16p^3+16p^2-4p)x^3-24px)P_{p-1}(x)^2
\\&\q+((16p^4-24p^2+5)x^4-(48p^2+4)x^2+24)P_p(x)P_{p-1}(x)\Big).
\end{align*}
\par{\it Proof.} Taking $y=-x$ in Theorem 2.2 and noting that $P_n(-x)=(-1)^nP_n(x)$ yields the result.
\v2
\par{\bf Corollary 2.5} {\sl Let $p\in\Bbb Z^+$ and $x\not=1$. Then}
\begin{align*}
&\sum_{n=0}^{p-1} (2n+1)^5 P_n(x)
 \\&=\f p{(1-x)^3}\Big(\big((2p+1)^4 (1-x)^2 - (32(2p+1)^2+48)(1-x) + 256\big) P_{p-1}(x) \\
&\q-\bigl((2p-1)^4 (1-x)^2 -(32(2p-1)^2+48)(1-x)+256\big)P_p(x)
\Big).
\end{align*}
 \vskip0.2cm
 \par{\it Proof.} Taking $y=1$ in Theorem 2.2 and noting that $P_n(1)=1$ gives the result.
\vskip0.2cm
\par{\bf Corollary 2.6} {\sl For $p=1,3,5,\ldots$ we have
$$\sum_{n=0}^{(p-1)/2}(4n+1)^5\f{\b{2n}n}{(-4)^n}
=(16p^5+32p^4-104p^3-120p^2+177p)(-1)^{\f{p-1}2}\f{\b{p-1}{(p-1)/2}}
{2^{p-1}}.$$
Hence, for any prime $p>3$ we have}
$$\sum_{n=0}^{(p-1)/2}(4n+1)^5\f{\b{2n}n}{(-4)^n}\e
2^{p-1}(32p^4-104p^3-120p^2+177p)+\f{59}4p^4B_{p-3}\mod {p^5}.$$
\par{\it Proof.} Putting $x=0$ in Corollary 2.5 and then applying
(1.2) and (2.6) yields the result.
\vskip0.2cm

\par{\bf Theorem 2.3} {\sl Let $p\in\Bbb Z^+$ and $x\not=\pm 1$. Then}
\begin{align*}&\sum_{n=0}^{p-1}(2n+1)^5P_n(x)^2
\\&=-\f{p^2}{15(x^2-1)^3}
\Big(\big((48p^2(p-1)^2+104p(p-1)+15)(x^2-1)^2
\\&\q+(256p(p-1)+400)(x^2-1)+512\big)
P_p(x)^2
\\&\q+\big((48p^2(p+1)^2+104p(p+1)+15)(x^2-1)^2
\\&\q+(256p(p+1)+400)(x^2-1)+512\big)
P_{p-1}(x)^2
\\&\q-2\big((48p^4-40p^2+7)(x^2-1)^2+(256p^2+144)(x^2-1)+512
\big)xP_p(x)P_{p-1}(x)\Big).
\end{align*}
\par{\it Proof.} Let $f_p(x)$ be the right-hand side of the identity in Theorem 2.3. Using the relation $(n+1)P_{n+1}(x)=(2n+1)xP_n(x)-nP_{n-1}(x)$ one may deduce that
$$f_{n+1}(x)-f_n(x)=(2n+1)^5P_n(x)^2.$$
Thus,
\begin{align*}\sum_{n=0}^{p-1}(2n+1)^5P_n(x)^2
=\sum_{n=0}^{p-1}(f_{n+1}(x)
-f_n(x))=f_p(x)-f_0(x)=f_p(x).\end{align*}
This proves the theorem.
\v2
\par{\bf Theorem 2.4} {\sl Let $p\in\Bbb Z^+$. Then}
$$(1-x)^4\sum_{n=0}^{p-1} (2n+1)^7 P_n(x)=pL_3(p,1-x)P_{p-1}(x)-pL_3(-p,1-x)P_p(x),$$
where
\begin{align*}L_3(p,t)&=(2p+1)^6t^3 - (72(2p+1)^4+400(2p+1)^2+256)t^2
 \\&\q+ (2304(2p+1)^2+6656)t - 18432.\end{align*}
 \vskip0.2cm
 \par{\it Proof.} Set
 $$ g(n)=nL_3(n,1-x)P_{n-1}(x)-nL_3(-n,1-x)P_n(x).$$
 Since $L_3(-n-1,t)=L_3(n,t)$ we see that
 \begin{align*} g(n+1)-g(n)&
 =(n+1)L_3(n+1,1-x)P_{n}(x)-(n+1)L_3(-n-1,1-x)P_{n+1}(x)
 \\&\q -nL_3(n,1-x)P_{n-1}(x)+nL_3(-n,1-x)P_n(x)
 \\&=(n+1)L_3(n+1,1-x)P_{n}(x)-L_3(n,1-x)((2n+1)xP_n(x)-nP_{n-1}(x))
 \\&\q -nL_3(n,1-x)P_{n-1}(x)+nL_3(-n,1-x)P_n(x)
 \\&=((n+1)L_3(n+1,1-x)-(2n+1)xL_3(n,1-x)+nL_3(-n,1-x))P_n(x)
 \\&=(1-x)^4(2n+1)^7P_n(x).
  \end{align*}
Thus,
$$(1-x)^4\sum_{n=0}^{p-1} (2n+1)^7 P_n(x)
=\sum_{n=0}^{p-1}(g(n+1)-g(n))=g(p)-g(0)=g(p).$$
This proves the theorem.
\v2

\par{\bf Theorem 2.5} {\sl Let $p\in\Bbb Z^+$. Then}
$$(1-x)^5\sum_{n=0}^{p-1} (2n+1)^9 P_n(x)=pL_4(p,1-x)P_{p-1}(x)-pL_4(-p,1-x)P_p(x),$$
where
\begin{align*}L_4(p,t)&= (2p+1)^8t^4-32(4(2p+1)^6 + 49(2p+1)^4 + 112(2p+1)^2 + 40)t^3 \\&\q+1536(6(2p+1)^4 + 66(2p+1)^2 + 89)t^2\\&\q-73728(4(2p+1)^2 + 17)t +2359296.\end{align*}
 \vskip0.2cm
 \par{\it Proof.} Set
 $$ g(n)=nL_4(n,1-x)P_{n-1}(x)-nL_4(-n,1-x)P_n(x).$$
 Since $L_4(-n-1,t)=L_4(n,t)$ we see that
 \begin{align*} g(n+1)-g(n)&
 =(n+1)L_4(n+1,1-x)P_{n}(x)-(n+1)L_4(-n-1,1-x)P_{n+1}(x)
 \\&\q -nL_4(n,1-x)P_{n-1}(x)+nL_4(-n,1-x)P_n(x)
 \\&=(n+1)L_4(n+1,1-x)P_{n}(x)-L_4(n,1-x)((2n+1)xP_n(x)-nP_{n-1}(x))
 \\&\q -nL_4(n,1-x)P_{n-1}(x)+nL_4(-n,1-x)P_n(x)
 \\&=((n+1)L_4(n+1,1-x)-(2n+1)xL_4(n,1-x)+nL_4(-n,1-x))P_n(x)
 \\&=(1-x)^5(2n+1)^9P_n(x).
  \end{align*}
Thus,
$$(1-x)^5\sum_{n=0}^{p-1} (2n+1)^9 P_n(x)
=\sum_{n=0}^{p-1}(g(n+1)-g(n))=g(p)-g(0)=g(p).$$
This proves the theorem.
\v2
\par{\bf Conjecture 2.1} {\sl Suppose that $m,p\in\Bbb Z^+$. Then there are integral polynomials $f_i(t)$ with degree $i$ $(i=0,1,\ldots,m-1)$ such that
$$(1-x)^{m+1}\sum_{n=0}^{p-1}(2n+1)^{2m+1}P_n(x)
=pL_m(p,1-x)P_{p-1}(x)-pL_m(-p,1-x)P_p(x),$$
where}
$$L_m(p,t)=(2p+1)^{2m}t^m+f_{m-1}((2p+1)^2)t^{m-1}
+\cdots+f_1((2p+1)^2)t+f_0.$$
\v2
\par For two numbers $b$ and $c$, following Z.W. Sun we define
$$T_n(b,c)=\sum_{k=0}^{[n/2]}\b n{2k}\b{2k}kb^{n-2k}c^k\ (n=0,1,2,\ldots).$$
By [20],
$$T_n(b,c)=(\sqrt{b^2-4c})^nP_n\Ls b{\sqrt{b^2-4c}}\q\t{for $b^2-4c\not=0$}.\eqno{(2.8)}$$
It is easy to show that (see [20, Lemma 2.5]) for $b,c\in\Bbb Z$ and any odd prime $p$,
$$T_p(b,c)\e b\mod p\q\t{and}\q T_{p-1}(b,c)\e \Ls{b^2-4c}p\mod p.\eqno{(2.9)}$$
\par From the above results for $P_n(x)$ one may easily deduce the identities for $T_n(b,c)$. For instance, we have the following result.
\v2
\par{\bf Theorem 2.6} {\sl Let $b$ and $c$ be two complex numbers with $bc(b^2-4c)\not=0$. Then}
\begin{align*}&\sum_{n=0}^{p-1}(2n+1)^3\f{T_n(b,c)^2}{(b^2-4c)^n}
\\&=\f{p^2}{12c^2(b^2-4c)^{p-1}}
\big(-((4p(p-1)-5)c+2b^2))T_p(b,c)^2
\\&\q-(4p(p+1)-5)c+2b^2))(b^2-4c)T_{p-1}(b,c)^2
+((8bcp^2-18bc+4b^3)T_p(b,c)T_{p-1}(b,c)\big)
\end{align*}
and
\begin{align*}\sum_{n=0}^{p-1}(2n+1)^3\f{T_n(b,c)^2}{(4c-b^2)^n}
&= (-1)^{p-1} \frac{p}{b^3(b^2-4c)^{p-1}} \big( 2p b (b^2-4c) T_{p-1}(b,c)^2 - 2p b T_p(b,c)^2 \\&\q+ ((4p^2-3)b^2+8c) T_p(b,c) T_{p-1}(b,c) \big).
\end{align*}
\par{\it Proof.} Taking $x=\f{b}{\sqrt{b^2-4c}}$ in Theorem 2.1(i) and Corollary 2.3 and then applying
 (2.8) yields the result.
\v2
\par{\bf Corollary 2.7} {\sl Let $p$ be a prime greater than $3$, $b,c\in\Bbb Z$ and $p\nmid bc(b^2-4c)$. Then
$$\sum_{n=0}^{p-1}(2n+1)^3\f{T_n(b,c)^2}{(b^2-4c)^n}
\e\begin{cases} -\f 53p^2\mod {p^3}&\t{if $\sls{b^2-4c}p=1$,}
\\-\f{2b^4-9b^2c+5c^2}{3c^2}p^2\mod{p^3}&\t{if $\sls{b^2-4c}p=-1$}
\end{cases}$$
and}
$$\sum_{n=0}^{p-1}(2n+1)^3\f{T_n(b,c)^2}{(4c-b^2)^n}
\e \f{-3b^2+8c}{b^2}\Ls{b^2-4c}pp\mod {p^2}.$$
\par{\it Proof.} This is immediate from Theorem 2.6 and (2.9).
\v2

\section*{3. Identities involving Ap\'ery-like numbers}
\par\q Let $r(n)$ and $s(n)$ be two sequences such that $r(n) \neq 0$ for $n = 0, 1, 2, \dots$. Suppose that $c \neq 0$ and $\{u_n\}$ is defined by
\begin{align*} u_{-1}=0,\ u_0=1\q\t{and}\q
r(n+1)u_{n+1}=s(n)u_n-cr(n)u_{n-1}\quad (n\ge 0).\notag\end{align*}
From  [17, (22.2)-(22.3)], for any positive integer $p$,
\begin{align*}\sum_{n=0}^{p-1} \bigl(r(n)x^2 - s(n)x + cr(n+1)\bigr)\f{u_n}{x^n} =\f{r(p)}{x^{p-1}}\bigl(cu_{p-1}-xu_p\bigr)+x^2r(0)  \tag{3.1}
\end{align*}
and
\begin{align*}
&\sum_{n=0}^{p-1} \bigl((n-3)r(n)x^2 - (n-2)s(n)x + c(n-1)r(n+1)\bigr)\frac{u_n}{x^n}\\
&=\frac{r(p)}{x^{p-1}} \bigl((p-2)cu_{p-1}-(p-3)xu_p\bigr)-3x^2r(0). \tag{3.2}
\end{align*}
Taking the derivative of $x$ on both sides of (3.2)
and then multiplying both sides by $x$ yields
\begin{align*}
&\sum_{n=0}^{p-1} \bigl((n-2)(n-3)r(n)x^2 - (n-1)(n-2)s(n)x + cn(n-1)r(n+1)\bigr)\frac{u_n}{x^n}
\\&= \frac{r(p)}{x^{p-1}} (p-2)\bigl((p-1)cu_{p-1}-(p-3)xu_p\bigr)+6x^2r(0). \tag{3.3}
\end{align*}
\par{\bf Theorem 3.1} {\sl
 Let $p\in\Bbb Z^+$, $a,b,d\in\Bbb R$, $d\not=0$ and $u_n=u_n(a,b,d^2)$. Then}
$$(d-a)\sum_{n=0}^{p-1}(2n+1)^3\frac{u_n}{d^n} =
\frac{4p^3}{d^{p-1}}\bigl(du_{p-1}-u_p\bigr) -(a-4b+3d)\sum_{n=0}^{p-1}(2n+1)\frac{u_n}{d^n}. $$

\par{\it Proof.}
Putting $r(n)=n^3$, $s(n)=(2n+1)\bigl(an(n+1)+b\bigr)$ and $c=d^2$ in (3.1) gives
$$\sum_{n=0}^{p-1}(n^3d^2-(2n+1)(an(n+1)+b)d+(n+1)^3d^2)
\f{u_n}{d^n}
=\f{p^3}{d^{p-1}}(d^2u_{p-1}-du_p).$$
That is,
$$\sum_{n=0}^{p-1}(2n+1)((n^2+n+1)d-a(n^2+n)-b)\f{u_n}{d^n}
=\f{p^3}{d^{p-1}}(du_{p-1}-u_p).$$
Hence,
$$\sum_{n=0}^{p-1}(2n+1)((2n+1)^2(d-a)+3d+a-4b)\f{u_n}{d^n}
=\f{4p^3}{d^{p-1}}(du_{p-1}-u_p),$$
which yields the result.
\vskip0.2cm

\par{\bf Corollary 3.1} {\sl For any positive integer $p$ we have}
$$2x(x+1)\sum_{n=0}^{p-1}(2n+1)V_n(x)=p^3(V_p(x)-V_{p-1}(x))$$
and
\begin{align*}
&\sum_{n=0}^{p-1}(2n+1)^3(-1)^nV_n(x)
\\&=2(-1)^{p-1}p^3(V_{p-1}(x)+V_p(x))
-(4x^2+4x+3)\sum_{n=0}^{p-1}(2n+1)(-1)^nV_n(x).\end{align*}
\par{\it Proof.} Putting $a=1$, $b=2x^2+2x+1$, $d=1$ and $u_n=V_n(x)$ in Theorem 3.1 gives the first identity, and taking $a=1$, $b=2x^2+2x+1$, $d=-1$ and $u_n=V_n(x)$ in Theorem 3.1 gives the second identity.

\vskip0.2cm
\par {\bf Theorem 3.2} {\sl
  Let $p\in\Bbb Z^+$, $a,b,d\in\Bbb R$, $d\not=0$ and $u_n=u_n(a,b,d^2)$. Then}
\begin{align*}
(a-d)^2\sum_{n=0}^{p-1}(2n+1)^5\frac{u_n}{d^n}
&=\frac{8p^3}{d^{p-1}}
\Bigl(\bigl( 2p(p+1)(d-a) - (9d+a-2b) \bigr)d u_{p-1}\\
&\q- \bigl( 2p(p-1)(d-a)-(9d+a-2b) \bigr) u_p
\Bigr)
\\&\quad + (41d^2 + (38a-88b)d + (a-4b)^2)
 \sum_{n=0}^{p-1} (2n+1) \frac{u_n}{d^n}.
\end{align*}

\par{\it Proof.}
Set
$$g(n)=\f{8n^3}{d^{n-1}}(2n(n+1)(d-a)
-(9d+a-2b))du_{n-1}-(2n(n-1)(d-a)-(9d+a-2b))u_n)$$
and
$$f(n)=\big((d-a)^2(2n+1)^5-(41d^2+(38a-88b)d+(a-4b)^2)
(2n+1)\big)\f{u_n}{d^n}.$$
Using (1.4) we see that
\begin{align*}&g(n+1)-g(n)
\\&=\f{8(n+1)^3}{d^{n}}\big(2(n+1)(n+2)(d-a)
-(9d+a-2b)\big)du_{n}
\\&\q-(2(n+1)n(d-a)-(9d+a-2b))u_{n+1})
\\&\q-\f{8n^3}{d^{n-1}}(2n(n+1)(d-a)
-(9d+a-2b))du_{n-1}-(2n(n-1)(d-a)-(9d+a-2b))u_n)
\\&=\f{8(n+1)^3}{d^{n}}(2(n+1)(n+2)(d-a)-(9d+a-2b))du_{n}
\\&\q-\f 8{d^n}(2n(n+1)(d-a)-(9d+a-2b))((2n+1)(an(n+1)+b)u_n
-n^3d^2u_{n-1})
\\&\q-\f{8n^3}{d^n}(2n(n+1)(d-a)-(9d+a-2b))d^2u_{n-1}
\\&\q+\f{8n^3}{d^n}(2n(n-1)(d-a)-(9d+a-2b))du_n
\\&=f(n).
\end{align*}
Thus,
$$\sum_{n=0}^{p-1}f(n)=\sum_{n=0}^{p-1}(g(n+1)-g(n))
=g(p)-g(0)=g(p).$$
\vskip0.2cm
\par{\bf Corollary 3.2} {\sl For any positive integer $p$ we have}
\begin{align*}&\sum_{n=0}^{p-1}(2n+1)^5(-1)^nV_n(x)
\\&=2p^3\big(4p(p+1)-(4x^2+4x+10))V_{p-1}(x)+(4p(p-1)-(4x^2+4x+10))V_p(x)\big)
\\&\q+(16x^4+32x^3+72x^2+56x+25)
\sum_{n=0}^{p-1}(2n+1)(-1)^nV_n(x).\end{align*}

\par{\it Proof.} Putting $a=1$, $b=2x^2+2x+1$, $d=-1$ and $u_n=V_n(x)$ in Theorem 3.2 gives the result.
\vskip0.2cm
\par{\bf Theorem 3.3} {\sl
  Let $p\in\Bbb Z^+$, $a,b,d\in\Bbb R$, $d\not=0$ and $u_n=u_n(a,b,d^2)$. Then
$$(d-a)^3\sum_{n=0}^{p-1}(2n+1)^7\frac{u_n}{d^n}
=\f{4p^3}{d^{p-1}}(\lambda(p)du_{p-1}-\lambda(-p)u_p)-
C\sum_{n=0}^{p-1}(2n+1)
\f{u_n}{d^n},$$
where
\begin{align*}\lambda(p)&=16(d-a)^2p^4+32(d-a)^2p^3
-4(d-a)(7a-4b+45d)p^2
\\&\q-4(d-a)(3a-4b+49d)p
+3a^2-12ab+16b^2+(310a-276b)d+727d^2\end{align*}
and
$$C=1595d^3+(2203a-3884b)d^2
+(361a^2-1416ab+1168b^2)d+(a-4b)^3.$$}
\par{\it Proof}. Set
$$g(n)=\f{4n^3}{d^{n-1}}(\lambda(n)du_{n-1}-\lambda(-n)u_n)
\ \t{and}\ f(n)=((d-a)^3(2n+1)^7+C(2n+1))\f{u_n}{d^n}.$$
It is easy to check that $\lambda(-n-1)=\lambda (n)$. Thus,
\begin{align*}&g(n+1)-g(n)\\&=\f{4(n+1)^3}{d^{n}}
(\lambda(n+1)du_{n}
-\lambda(-n-1)u_{n+1})
-\f{4n^3}{d^{n-1}}(\lambda(n)du_{n-1}-\lambda(-n)u_n)
\\&=\f{4(n+1)^3}{d^{n-1}}\lambda(n+1)u_n
-\f{4\lambda(-n-1)}{d^n}\big((2n+1)(an(n+1)+b)u_n
-d^2n^3u_{n-1}\big)
\\&\q-\f{4n^3}{d^{n-1}}\lambda(n)du_{n-1}
+\f{4n^3}{d^{n-1}}\lambda(-n)u_n
\\&=4\big((n+1)^3\lambda(n+1)d
-\lambda(n)(2n+1)(an(n+1)+b)+n^3\lambda(-n)d\big)\f{u_n}{d^n}
\end{align*}
One can check that
\begin{align*}&4((n+1)^3\lambda(n+1)d
-\lambda(n)(2n+1)(an(n+1)+b)+n^3\lambda(-n)d)
\\&=(d-a)^3(2n+1)^7
+C(2n+1).\end{align*}
Thus, $g(n+1)-g(n)=f(n)$ and so
$$\sum_{n=0}^{p-1}f(n)=\sum_{n=0}^{p-1}(g(n+1)-g(n))
=g(p)-g(0)=g(p),$$
which yields the result.
\vskip0.2cm
\par{\bf Corollary 3.3} {\sl For any positive integer $p$ we have}
\begin{align*}&\sum_{n=0}^{p-1}(2n+1)^7(-1)^nV_n(x)
\\&=2p^3\big(\lambda_0(p) V_{p-1}(x) + \lambda_0(-p) V_p(x)\big)
\\&\q-(64x^6+192x^5+848x^4+1376x^3+1884x^2+1228x+427)
\sum_{n=0}^{p-1}(2n+1)(-1)^nV_n(x),\end{align*}
where
$$\lambda_0(p)=16p^4+32p^3-(16x^2+16x+84)p^2-(16x^2+16x+100)p
+16x^4+32x^3+164x^2+148x+175.$$
\vskip0.2cm
\par{\it Proof.} Putting $a=1$, $b=2x^2+2x+1$, $d=-1$ and $u_n=V_n(x)$ in Theorem 3.3 gives the result.
\vskip0.2cm

\section*{4. Supercongruences involving Ap\'ery-like numbers}
\par\q\  Using the identities in Section 3 and known congruences for Ap\'ery-like numbers, in this section we deduce new supercongruences involving the first kind of Ap\'ery-like numbers.
\par {\bf Lemma 4.1 ([7])} {\sl
Let $p > 3$ be a prime. Then}
$$A_{p-1} \equiv 1 + \frac{2}{3} p^3 B_{p-3} \pmod{p^4}
\q\t{and}\q
A_p \equiv 5 - \frac{14}{3} p^3 B_{p-3} \pmod{p^4}.$$

\par{\bf Lemma 4.2 ([19])} {\sl
Let $p > 3$ be a prime. Then}
\[\sum_{n=0}^{p-1} (2n+1)A_n \equiv p+\frac{7}{6}p^4B_{p-3} \pmod{p^5}.\]

\par{\bf Theorem 4.1} {\sl
Let $p > 3$ be a prime. Then}
$$\sum_{n=0}^{p-1} (2n+1)^5 A_n \equiv p + 2p^3 + p^4\Big( -6+\frac{7}{6}B_{p-3}\Big)\pmod{p^5}$$
and
$$\sum_{n=0}^{p-1} (2n+1)^7 A_n \equiv 8p + 19p^3  +p^4 \Big( -30+\frac{28}{3}B_{p-3}\Big)  \pmod{p^5}.$$

\par{\it Proof.}
Putting $a=17$, $b=5$, $d=1$ and $u_n=A_n$ in Theorem 3.2 and then applying Lemmas 4.1 and 4.2 we derive

\begin{align*}
(-16)^2\sum_{n=0}^{p-1} (2n+1)^5 A_n
&= 8p^3\big((2p(p+1)(-16)-16)A_{p-1}-(2p(p-1)(-16)-16)A_p\big)
\\&\q+(41+38\cdot 17-88\cdot 5+(17-20)^2)\sum_{n=0}^{p-1}(2n+1)A_n
\\&\e 128p^3\big(-2p-1
+5(-2p+1)\big)
+256\Big(p+\f 76p^4B_{p-3}\Big)
\\&=512p^3+256p+256p^4\Big(-6+\f 76B_{p-3}\Big)\mod {p^5},
\end{align*}
which yields the first result.
\par Taking $a=17$, $b=5$, $d=1$ and $u_n=A_n$ in Theorem 3.3 and then applying Lemmas 4.1 and 4.2, we deduce that
\begin{align*}
\sum_{n=0}^{p-1} (2n+1)^7 A_n
&=-\frac{p^3}{4}(16p^4+32p^3+36p^2+20p+19) A_{p-1} \\ &\q+\frac{p^3}{4}(16p^4-32p^3+36p^2-20p+19) A_p+ 8\sum_{n=0}^{p-1} (2n+1) A_n
\\&\e -\frac{p^3}{4}(20p+19)+\frac{p^3}{4}(-20p+19)\cdot 5+8\Big(p+\f 76p^4B_{p-3}\Big)
\\&=8p + 19p^3  + \Big( \frac{28}{3}B_{p-3} - 30 \Big) p^4  \pmod{p^5}.
\end{align*}
The proof is now complete.
\v2
\par{\bf Remark 4.1}
Let $p > 3$ be a prime. In [13], Mao and Ni showed that
$$\sum_{n=0}^{p-1}(2n+1)^3A_n\e p^3-\f 43p^6B_{p-3}\mod {p^7}.$$
In [3], Guo and Zeng confirmed the following conjecture due to Z.W. Sun:
$$\sum_{n=0}^{p-1} (2n+1)(-1)^n A_n \equiv p\Big( \frac{p}{3} \Big) \pmod{p^3}.$$
In [21], Z.W. Sun proved that
\begin{align*}
&\sum_{n=0}^{p-1} (2n+1)^3(-1)^n A_n \equiv -\frac{p}{3}\Big( \frac{p}{3} \Big) \pmod{p^3},
\\&\sum_{n=0}^{p-1} (2n+1)^5(-1)^n A_n \equiv -\frac{13}{27}p\Ls p3 \pmod{p^3},
\\&
\sum_{n=0}^{p-1} (2n+1)^7(-1)^n A_n \equiv\frac{5}{9}p\left( \frac{p}{3} \right) \pmod{p^3}.
\end{align*}

\par{\bf Lemma 4.3 ([5])} {\sl
Let $p > 3$ be a prime. Then}
\[
\sum_{n=0}^{p-1} (2n+1)\frac{b_n}{(-9)^n} \equiv p\left( \frac{p}{3} \right) \pmod{p^3}.
\]

\par{\bf Lemma 4.4 ([8])} {\sl
Let $p > 3$ be a prime. Then}
\begin{align*}
&b_{p-1} \equiv 81^{p-1} - \frac{2}{27}p^3 B_{p-3} \pmod{p^4}, \\&b_p \equiv -3 - 6p^3 B_{p-3} \pmod{p^4}.
\end{align*}

\par{\bf Theorem 4.2} {\sl
Let $p > 3$ be a prime. Then
$$\sum_{n=0}^{p-1} (2n+1)^5\frac{b_n}{(-9)^n} \equiv 841p\left( \frac{p}{3} \right)  \pmod{p^3}$$
and
$$\sum_{n=0}^{p-1} (2n+1)^7\frac{b_n}{(-9)^n} \equiv  -181763\, p \left(\frac{p}{3}\right) \pmod{p^3}.$$}
\par{\it Proof.} Putting $a=-7$, $b=-3$, $d=-9$ and $u_n=b_n$ in Theorem 3.2 and then applying Lemmas 4.3 and 4.4 gives
\begin{align*}
4\sum_{n=0}^{p-1} (2n+1)^5 \frac{b_n}{(-9)^n}
&= \f {8p^3}{9^{p-1}}\big((-9)(-4p(p+1)-82)b_{p-1}+(4p(p-1)-82)b_p\big)
\\&\q+(41\cdot81+38\cdot (-7)-88\cdot (-3)+(-7+12)^2)\sum_{n=0}^{p-1}(2n+1)\frac{b_n}{(-9)^n}\\
&\e 3364p \Big( \frac{p}{3} \Big) \pmod{p^3},
\end{align*}
which yields the first result.
\par Taking $a=-7$, $b=-3$, $d=-9$ and $u_n=b_n$ in Theorem 3.3 and then applying Lemmas 4.3 and 4.4 deduces that
\begin{align*}
&\sum_{n=0}^{p-1} (2n+1)^7 \frac{b_n}{(-9)^n}
\\&=\frac{p^3}{(-9)^p}\big((-2592p^4 - 5184p^3 + 143208p^2 + 145800p - 2875662)\,b_{p-1}
 \\&\q+ (-288p^4 + 576p^3 + 15912p^2 - 16200p - 319518)\,b_p \big)
- 181763 \sum_{n=0}^{p-1} (2n+1)\frac{b_n}{(-9)^n}
\\&\equiv  -181763\, p \Big(\frac{p}{3}\Big) \pmod{p^3}.
\end{align*}
The proof is now complete.
\vskip0.2cm
\par{\bf Lemma 4.5 ([5])} {\sl
Let $p > 3$ be a prime. Then}
\[
\sum_{n=0}^{p-1} (2n+1)\frac{T_n}{4^n} \equiv p+\frac{7}{6}p^4B_{p-3} \pmod{p^5}.
\]

\par{\bf Lemma 4.6 ([6,12])} {\sl
Let $p > 3$ be a prime. Then}
$$T_{p-1} \equiv 16^{p-1} + \frac{p^3}{4} B_{p-3} \pmod{p^4}
\q\t{and}\q
T_p \equiv 4 - p^3 B_{p-3} \pmod{p^4}.$$

\par{\bf Theorem 4.3} {\sl
Let $p > 3$ be a prime and $q_p(2)=(2^{p-1}-1)/p$. Then
$$\sum_{n=0}^{p-1} (2n+1)^5\frac{T_n}{4^n} \equiv  17p  + p^4\Big(\frac{119}{6}B_{p-3}-80q_p(2)-16\Big) \pmod{p^5}$$
and
\begin{align*}
\sum_{n=0}^{p-1} (2n+1)^7\frac{T_n}{4^n}
&\equiv  561p-432p^4
+ p^4\Big(\frac{1309}{2}B_{p-3}-2776q_p(2)\Big)
 \pmod{p^5}.
\end{align*}}

\par{\it Proof.}
Putting $a=12$, $b=4$, $d=4$ and $u_n=T_n$ in Theorem 3.2 and then applying Lemmas 4.5-4.6 and the fact that $4^{p-1}\e 1+2pq_p(2)\mod{p^2}$ and $4^{-(p-1)}\e 1-2pq_p(2)\mod {p^2}$  
yields
\begin{align*}
&(4-12)^2\sum_{n=0}^{p-1} (2n+1)^5\frac{T_n}{4^n}
\\&= (41\cdot 4^2+(38\cdot 12-88\cdot 4)\cdot 4+(12-4\cdot 4)^2)
\sum_{n=0}^{p-1} (2n+1)\frac{T_n}{4^n}
\\&\q+ \frac{8p^3}{4^{p-1}} \big((2p(p+1)(-8)-40)\cdot 4T_{p-1}-(2p(p-1)(-8)-40)T_p\big)
 \\
&\equiv 64\cdot 17\Big(p+\f 76p^4B_{p-3}\Big)
+\f{8p^3}{4^{p-1}}\big((-16p-40)\cdot 4\cdot 16^{p-1}-(16p-40)\cdot 4\big)
\\&\e 64\cdot 17\Big(p+\f 76p^4B_{p-3}\Big)+64p^3\big(-16p-80pq_p(2)\big)
\pmod{p^5},
\end{align*}
which yields the first result.
\par Taking $a=12$, $b=4$, $d=4$ and $u_n=T_n$ in Theorem 3.3 and then applying Lemmas 4.5-4.6 yields
\[
\begin{aligned}
\sum_{n=0}^{p-1} (2n+1)^7\frac{T_n}{4^n}
&= \frac{p^3}{4^p} \big(
\big(-128p^4 - 256p^3 - 992p^2 - 864p - 2776\big)T_{p-1} \\
&\quad + \big(32p^4 - 64p^3 + 248p^2 - 216p + 694\big)T_p\big)
+ 561 \sum_{n=0}^{p-1} (2n+1)\frac{T_n}{4^n} \\
&\e   561p + \frac{1309}{2}\,p^4 B_{p-3}
+ \frac{p^3}{4^p}\big(2776\bigl(1-16^{p-1}\bigr) - 864p\bigl(1+16^{p-1}\bigr)\big)\\
&\e 561p-432p^4
+ p^4\Big(\frac{1309}{2}B_{p-3}-2776q_p(2)\Big)
\pmod{p^5}.
\end{aligned}
\]
The proof is now complete.
\vskip0.2cm

\par{\bf Lemma 4.7 ([4])} {\sl
Let $p > 3$ be a prime. Then}
\[
\sum_{n=0}^{p-1} (2n+1)\frac{T_n}{(-4)^n} \equiv (-1)^{\frac{p-1}{2}}p +p^3E_{p-3}\pmod{p^4}.
\]

\par{\bf Theorem 4.4} {\sl
Let $p > 3$ be a prime. Then
$$\sum_{n=0}^{p-1} (2n+1)^5\frac{T_n}{(-4)^n} \equiv (-1)^{\frac{p-1}{2}}p + p^3\big(E_{p-3} - 8\big) \pmod{p^4}$$
and
$$\sum_{n=0}^{p-1} (2n+1)^7\frac{T_n}{(-4)^n}  \e 15(-1)^{\f {p-1}{2}} p + p^3 \bigl(10 + 15E_{p-3}\bigr) \pmod{p^4}.$$}
\par{\it Proof.}
Putting $a=12$, $b=4$, $d=-4$ and $u_n=T_n$ in Theorem 3.2 and then applying Lemmas 4.6-4.7 yields
\begin{align*}
\sum_{n=0}^{p-1} (2n+1)^5 \frac{T_n}{(-4)^n}
&= \frac{p^3}{(-4)^{p-1}} \big( 4(p^2 + p - 1) T_{p-1} + (p^2 - p - 1) T_p \big) + \sum_{n=0}^{p-1} (2n+1) \frac{T_n}{(-4)^n}
 \\&\e\frac{p^3}{(-4)^{p-1}} (-4\cdot16^{p-1}-4)+(-1)^{\frac{p-1}{2}}p +p^3E_{p-3}
 \\&\e (-1)^{\frac{p-1}{2}} p + p^3 \big( E_{p-3} - 8 \big) \pmod{p^4},
\end{align*}
which yields the first result.
\par Taking $a=12$, $b=4$, $d=-4$ and $u_n=T_n$ in Theorem 3.3 and then applying Lemmas 4.6 and 4.7, we deduce that
\begin{align*}
\sum_{n=0}^{p-1} (2n+1)^7 \frac{T_n}{(-4)^n}
&=\frac{p^3}{(-4)^p} \big(\big(-64p^4 - 128p^3 + 112p^2 + 176p - 20\big)T_{p-1}\\
&\q+\big(-16p^4 + 32p^3 + 28p^2 - 44p - 5\big)T_p \big)+ 15 \sum_{n=0}^{p-1} (2n+1)\f{T_n}{(-4)^n}\\
&\e\frac{p^3}{(-4)^p} \big( -20\cdot16^{p-1}-20\big)+15(-1)^{\frac{p-1}{2}}p +15p^3E_{p-3}\\
&\e 15(-1)^{\frac{p-1}{2}} p + p^3 \bigl(10 + 15E_{p-3}\bigr)  \pmod{p^4}.
\end{align*}
The proof is now complete.
\vskip0.2cm
\par{\bf Lemma 4.8 ([10])} {\sl
Let $p > 3$ be a prime. Then}
\[
\sum_{n=0}^{p-1} (2n+1)\frac{D_n}{8^n} \equiv  p + \frac{35}{24} p^4 B_{p-3} \pmod{p^5}.
\]

\par{\bf Lemma 4.9 ([12,23])}
{\sl Let $p > 3$ be a prime. Then}
\begin{align*}
&D_{p-1} \equiv 64^{p-1} - \frac{p^3}{6} B_{p-3} \pmod{p^4},
\\&D_p \equiv 4 + \frac{16}{3} p^3 B_{p-3} \pmod{p^4}.
\end{align*}

\par{\bf Theorem 4.5} {\sl
Let $p > 3$ be a prime and $q_p(2)=(2^{p-1}-1)/p$. Then
$$\sum_{n=0}^{p-1} (2n+1)^5\frac{D_n}{8^n} \equiv 721p -592p^3+p^4\Big(\frac{25235}{24} B_{p-3}-96-5328q_p(2)\Big) \pmod{p^5}$$
and
\begin{align*}
&\sum_{n=0}^{p-1} (2n+1)^7 \frac{D_n}{8^n} \\&\equiv  152153p   - 125144p^3 + p^4\Big(\frac{5325355}{24} B_{p-3}-1126296q_p(2) - 19488 \Big)\pmod{p^5}.
\end{align*}}
\par{\it Proof.}
Putting $a=10$, $b=4$, $d=8$ and $u_n=D_n$ in Theorem 3.2 and then applying Lemmas 4.8-4.9 and the fact that $8^{p-1}=(1+pq_p(2))^3\e 1+3pq_p(2)\mod {p^2}$ and so
$8^{-(p-1)}\e 1-3pq_p(2)\mod {p^2}$
yields
\begin{align*}
&\sum_{n=0}^{p-1} (2n+1)^5\frac{D_n}{8^n}
\\&= \frac{2p^3}{8^{p-1}} \big(( -32p^2 - 32p - 592 ) D_{p-1} +(4p^2 - 4p + 74)D_p \big)+ 721\sum_{n=0}^{p-1} (2n+1)\frac{D_n}{8^n} \\
&\e \frac{2p^3}{8^{p-1}} \big(64^{p-1}(- 32p - 592)+4(- 4p + 74 )\big)+721\Big( p + \frac{35}{24} p^4 B_{p-3}\Big)\\
&\e 721p -592p^3+p^4\Big(\frac{25235}{24} B_{p-3}-96-5328q_p(2)\Big) \pmod{p^5}.
\end{align*}
\par Taking $a=10$, $b=4$, $d=8$ and $u_n=D_n$ in Theorem 3.3 and then applying Lemmas 4.8 and 4.9, we deduce that
\begin{align*}
&\sum_{n=0}^{p-1} (2n+1)^7 \frac{D_n}{8^n}
\\&=\frac{p^3}{8^p} \big(
\big(-2048p^4 - 4096p^3 - 105984p^2 - 103936p - 2002304\big)D_{p-1} \\
&\q+ \big(256p^4 - 512p^3 + 13248p^2 - 12992p + 250288\big)D_p
\big)+ 152153 \sum_{n=0}^{p-1} (2n+1)\frac{D_n}{8^n}\\
&\e \frac{p^3}{8^p} \big((- 103936p - 2002304)64^{p-1}+4(- 12992p + 250288)\big) + 152153\Big( p + \frac{35}{24} p^4 B_{p-3}\Big)\\
&\e 152153p   - 125144p^3 + p^4\Big(\frac{5325355}{24} B_{p-3}-1126296q_p(2) - 19488 \Big)\pmod{p^5}.
\end{align*}
The proof is now complete.
\vskip0.2cm
\par{\bf Lemma 4.10 ([16,22])} {\sl
Let $p > 3$ be a prime. Then}
\begin{align*}
&\sum_{n=0}^{p-1} (2n+1)\frac{V_n}{(-16)^n} \equiv (-1)^{\frac{p-1}{2}} p + 3p^3 E_{p-3} \pmod{p^4},\\
&\sum_{n=0}^{p-1} (2n+1)\frac{V_n^{(3)}}{(-27)^n} \equiv (-1)^{\left\lbrack \frac{p}{3} \right\rbrack} p + 7p^3 U_{p-3} \pmod{p^4}, \\
&\sum_{n=0}^{p-1} (2n+1)\frac{V_n^{(4)}}{(-64)^n} \equiv (-1)^{\left\lbrack \frac{p}{4} \right\rbrack} p + 13p^3 s_{p-3} \pmod{p^4}, \\
&\sum_{n=0}^{p-1} (2n+1)\frac{V_n^{(6)}}{(-432)^n} \equiv (-1)^{\frac{p-1}{2}} p + \frac{155}{9} p^3 E_{p-3} \pmod{p^4},
\end{align*}
where $\{s_n\}$ is given by
$s_0=1$ and $s_n=1-\sum_{k=0}^{n-1} \b
nk2^{2n-1-2k}s_k\ (n\ge 1).$
\v2
\par{\bf Lemma 4.11 ([9,11])} {\sl
Let $p > 3$ be a prime. Then}
\begin{align*}
&V_p \equiv 8 + 40p^3 B_{p-3} \pmod{p^4}, \\
&V_p^{(3)} \equiv 15 + 132p^3 B_{p-3} \pmod{p^4}, \\
&V_p^{(4)} \equiv 40 + 704p^3 B_{p-3} \pmod{p^4}, \\
&V_p^{(6)} \equiv 312 + 16120p^3 B_{p-3} \pmod{p^4}, \\
&V_{p-1} \equiv 256^{p-1} - \frac{3}{2}p^3 B_{p-3} \pmod{p^4}, \\
&V_{p-1}^{(3)} \equiv 729^{p-1} - \frac{92}{27}p^3 B_{p-3} \pmod{p^4}, \\
&V_{p-1}^{(4)} \equiv 4096^{p-1} - \frac{17}{2}p^3 B_{p-3} \pmod{p^4}, \\
&V_{p-1}^{(6)} \equiv 1866624^{p-1} - \frac{1705}{54}p^3 B_{p-3} \pmod{p^4}.
\end{align*}

\par{\bf Lemma 4.12 ([18, Theorem 3.6])} {\sl
Let $p > 3$ be a prime, $x\in\Bbb Z_p$, and let  $\langle x{\rangle }_{p} $ be given by $ x \equiv  \langle x{\rangle }_{p}\left( {\;\operatorname{mod}\;p}\right)$ and $ \langle x{\rangle }_{p} \in  \{ 0,1,\ldots ,p - 1\} $. Set $ x' = ( {x-\langle x{\rangle }_{p}})/p $. Then}
$$\sum_{n=0}^{p-1} (2n+1) (-1)^n V_n(x)
\e (-1)^{\langle x\rangle_p}p
+ p^3\bigl(x'(x'+1)+1\bigr)E_{p-3}(-x) \pmod{p^4}.$$
\par{\bf Theorem 4.6} {\sl
Let $p > 3$ be a prime, $x\in\Bbb Z_p$ and let $\langle x\rangle_p$ and $x'$ be as in Lemma 4.12. Then
\begin{align*}
\sum_{n=0}^{p-1} (2n+1)^5 (-1)^n V_n(x)
&\e (-1)^{\langle x \rangle_p}  \big(16x^4 + 32x^3 + 72x^2 + 56x + 25\big)p  \\
&\q + (x'^2+x'+1) \big(-16x^2-16x-40\\
&\q + \big(16x^4 + 32x^3 + 72x^2 + 56x + 25\big) E_{p-3}(-x) \big)p^3 \pmod{p^4}  \end{align*}
and
\begin{align*}
&\sum_{n=0}^{p-1}(2n+1)^7(-1)^n V_n(x)
 \\&\equiv -(-1)^{\langle x \rangle_p} ( 64x^6+192x^5+848x^4+1376x^3+1884x^2+1228x+427 )p
  \\&\q+ (x'^2+x'+1) \big( 64x^4+128x^3+656x^2+592x+700 \\ &\q - ( 64x^6+192x^5+848x^4+1376x^3+1884x^2+1228x+427) E_{p-3}(-x) \big) p^3 \pmod{p^4}.
\end{align*}}
\par{\it Proof.} Since $x\in\Bbb Z_p$ we see that
\begin{align*}\b xp&=\b{px'+\xp}p=\f{(px'+\xp)(px'+\xp-1)\cdots(px'+1)px'}{p!}
\\&\q\times(p(x'-1)+p-1)\cdots(p(x'-1)+\xp+1)
\\&\e x'\pmod p\end{align*}
and so $$\b xp\b{-1-x}p\e x'\cdot\f{-1-x-(p-1-\xp)}p
=-x'(x'+1)\pmod p.$$
Note that $\b pk=\f{p(p-1)\cdots(p-k+1)}{k!}\e 0\pmod p$
for $k=1,2,\ldots,p-1$. From the above we derive
\begin{align*} V_p(x)&=1-\b{2p}p\b
xp\b{-1-x}p
+\sum_{k=1}^{p-1}\b pk\b {p+k}k(-1)^k\b xk\b{-1-x}k
\\&\e 1-2\b{2p-1}{p-1}\b
xp\b{-1-x}p\e 1+2x'(x'+1)\pmod p.\end{align*}
On the other hand,
\begin{align*} V_{p-1}(x)&=\sum_{k=0}^{p-1}\b{p-1}k\b{p-1+k}k(-1)^k\b
xk\b{-1-x}k
\\&=1+\sum_{k=1}^{p-1}\f{p(p-k)(p^2-1^2)(p^2-2^2)\cdots(p^2-(k-1)^2)}{k!^2}
(-1)^k\b xk\b{-1-x}k
\\&\e 1\pmod p.\end{align*}
\par From the above, Corollaries 3.2-3.3 and Lemma 4.12 we deduce that
\begin{align*}&\sum_{n=0}^{p-1}(2n+1)^5(-1)^nV_n(x)
\\&\e 2p^3\big(-(4x^2+4x+10)-(4x^2+4x+10)(1+2x'(x'+1))\big)
\\&\q+(16x^4+32x^3+72x^2+56x+25)
\big((-1)^{\langle x\rangle_p}p
+ p^3(x'(x'+1)+1)E_{p-3}(-x)\big)
\mod{p^4}\end{align*}
and
 \begin{align*}\sum_{n=0}^{p-1}(2n+1)^7(-1)^nV_n(x)
&\e 2p^3(16x^4+32x^3+164x^2+148x+175)(1+1+2x'(x'+1))
\\&\q-(64x^6+192x^5+848x^4+1376x^3+1884x^2+1228x+427)
\\&\q\times\big((-1)^{\langle x\rangle_p}p
+ p^3(x'(x'+1)+1)E_{p-3}(-x)\big)
\mod {p^4},\end{align*}
 which yields the results.
 \vskip0.2cm
\par {\bf Remark 4.2} Let $p$ be an odd prime and $x\in\Bbb Z_p$. In [18], Sun established the congruence for $\sum_{n=0}^{p-1}(2n+1)V_n(x)$ modulo $p^5$. In [14], Mao and Yang obtained the congruences for
$$\sum_{n=0}^{p-1}(2n+1)^3V_n(x),\q
\sum_{n=0}^{p-1}(2n+1)^3(-1)^nV_n(x),\q
\sum_{n=0}^{p-1}(2n+1)^5V_n(x)$$
modulo $p^4$.

 \vskip0.2cm
\par{\bf Corollary 4.1} {\sl Let $p > 3$ be a prime. Then}
\begin{align*}
&\sum_{n=0}^{p-1}(2n+1)^5\frac{V_n}{(-16)^n}
\equiv 12(-1)^{\frac{p-1}{2}}p + p^3(-27+36E_{p-3}) \pmod{p^4}, \\
&\sum_{n=0}^{p-1}(2n+1)^5\frac{V_n^{(3)}}{(-27)^n}
\equiv \frac{1081}{81}(-1)^{\left\lbrack \frac{p}{3} \right\rbrack}p
+ p^3\Big(-\frac{2296}{81}+\frac{7567}{81}U_{p-3}\Big) \pmod{p^4}, \\
&\sum_{n=0}^{p-1}(2n+1)^5\frac{V_n^{(4)}}{(-64)^n}
\equiv \frac{241}{16}(-1)^{\left\lbrack \frac{p}{4} \right\rbrack}p
+ p^3\Big(-\frac{481}{16}+\frac{3133}{16}s_{p-3}\Big) \pmod{p^4}, \\
&\sum_{n=0}^{p-1}(2n+1)^5\frac{V_n^{(6)}}{(-432)^n}
\equiv \frac{1420}{81}(-1)^{\frac{p-1}{2}}p
+ p^3\Big(-\frac{2635}{81}+\frac{220100}{729}E_{p-3}\Big) \pmod{p^4}
\end{align*}
and
\begin{align*}
&\sum_{n=0}^{p-1}(2n+1)^7\frac{V_n}{(-16)^n} \equiv -160(-1)^{\frac{p-1}{2}}p + 417p^3 - 480p^3E_{p-3} \pmod{p^4},\\
&\sum_{n=0}^{p-1}(2n+1)^7\frac{V_n^{(3)}}{(-27)^n} \equiv -\frac{135451}{729}(-1)^{[\f p3]}p + \frac{324100}{729}p^3 - \frac{948157}{729}p^3U_{p-3} \pmod{p^4},\\
&\sum_{n=0}^{p-1}(2n+1)^7\frac{V_n^{(4)}}{(-64)^n} \equiv -\frac{14041}{64}(-1)^{[\f p4]}p + \frac{30745}{64}p^3 - \frac{182533}{64}p^3 s_{p-3} \pmod{p^4},\\
&\sum_{n=0}^{p-1}(2n+1)^7\frac{V_n^{(6)}}{(-432)^n} \equiv -\frac{196048}{729}(-1)^{\frac{p-1}{2}}p + \frac{388585}{729}p^3 - \frac{30387440}{6561}p^3E_{p-3} \pmod{p^4}.
\end{align*}
\par{\it Proof.} Taking $ x =  - \f{1}{2}, - \f{1}{3}, - \f{1}{4}, - \f{1}{6} $ in  Theorem 4.6 and then applying Lemmas 4.10-4.11 deduces the results.

\end{document}